# Failure localization in time critical market applications

Mikhail Davidson
Carana Corp., Moscow, Russia
SkolTech, Moscow, Russia
MSU, Moscow, Russia
mdavidson@carana-corp.com

Gleb Labutin
Markets Development Department
JSC "System Operator of United Power System"
Moscow, Russia
labutin@so-ups.ru

*Abstract*—In time critical market applications such as for example scheduling and price computation for the balancing market, failure of the algorithm in finding a solution would result in cancelation of the session and respective financial consequences for the market participants. In the paper we propose a regularization procedure that could help an operator to spot the problematic location and find out the reasons that caused the algorithm to breakdown and make necessary corrections. The approach has proved to be useful in the operation of the Balancing Market in the Russian Federation.

*Keywords—Energy market; ACOPF; load flow; loss function; marginal loss coefficients*

I. INTRODUCTION

Russian energy market consists of the Day Ahead and Balancing markets. The former computes the traded volumes and the corresponding locational marginal prices for every hour of the next day, while the latter is used to trade the deviations from the day-ahead volumes and also to produce generation plant schedules for the period until the end of the current day with the hourly breakdown. The market operation covers more than 95% of the Russian energy system except for some isolated or thinly connected regions [1]. The peak load of the area under market operation exceeds 150 GWt. The system is highly geographically expanded with significant inter area power transfers. This requires a more accurate modeling in the scheduling and market procedures. To ensure a reasonable level of accuracy the procedures are based on ACOPF and the grid model consists of around 9000 buses. According to the data published by the System Operator of United Power System of Russia the deviation of the actual dispatch of steam power stations from the schedule is within 1% [2].

Balancing Market runs every one hour ahead of real time. This time is spent to get the updated information on the system conditions, including the topology, system constraints, load forecast, generation resources and their bids, and the market scheduling run itself. Hence, if the scheduling procedure fails for some reason the time window to analyze the problem and take corrective measures is rather narrow. The failure means that the true market based schedules are not issued and the system continues to operate based on the latest previous schedules with the dispatcher manual corrections according to situation. This inevitably increases side payments, decreases transparency and every effort is made to avoid such occurrences. But if it nevertheless happens it is highly important to localize and handle the problem in order to minimize the potential negative effects.

Since ACOPF is a nonlinear nonconvex problem there could be numerous reasons leading to abnormal termination of the algorithm. Feasible solution may not exist for a particular nodal load profile, and even for simple examples feasible domain turns out to be nonconvex [3]. In the energy system under stressed conditions a feasible solution may exist but close to the boundary that could also result in an unstable behavior of the numerical methods.

In the literature there are a number of approaches that can help handling the problems related to infeasibility and poor convergence (see e.g., [4, 5] for a survey of solution techniques of ACOPF). The algorithm of "in the process" correction of the problem data is proposed in [6] that helps to keep solvability of the problem.

Another issue, however, is that in the case of the convergence problems there is often no a clear indication of a problematic location in the grid that caused the breakdown. The focus of this paper is to develop a mechanism that would provide certain guidance to the operating personnel as to identification or localization of the source of the problem in the energy system which is important in the context of time critical application.

To simplify the analysis we restrict ourselves to a model problem of finding a solution to load flow equations for a given vector of net nodal injections. Most of the problems causing divergence of an ACOPF solver would arise when solving load flow problem as well.

The approach is based on minimizing squared norm of the residual of the load flow equations regularized by the active power loss function multiplied by a regularization parameter. In the case of convergence problems adding this regularizing term with certain parameter value ensures that resulting goal function is convex and hence facilitates the convergence. Additionally, infeasibility is usually caused by a failure of the transmission system to provide inflow/outflow of a required amount of power to/from a particular node or area, thus such problems often have a local nature. We show that the stationary points of such a regularized function possess the property that the residual nodal load error tend to attain its maximum values in the vicinity of a problem location thus indicating the reason of the divergence of the solution procedure. Preliminary experience demonstrates efficiency of the approach in the Russian energy system Balancing Market computations.

## II. NOTATION, DEFINITIONS, AND PRELIMINARY RESULTS

Consider the system of load flow equations

$$F(x,S) = 0 \quad (1)$$

where

$$x = (u,v) \in R^{2n}, \quad S = (P,Q) \in R^{2n}$$

are respectively the vectors of the state variables and active and reactive components of nodal power injections with

$$U = u + iv \in Z^n$$

being the vector of complex nodal voltages and $n$ the number of load and generation buses in the system (the swing bus is considered to have number 0).

Define active power loss function as the sum of active power losses over all lines $(r,k)$ in the system:

$$L(x) = L(u,v) = \sum_{(r,k)} |U_r - U_k|^2 R_{rk}/(R_{rk}^2 + X_{rk}^2)$$

where $R_{rk}$ and $X_{rk}$ are line resistance and reactance respectively.

*Property 1.* If all line resistances $R_{rk}$ are strictly positive function $L(x)$ is strongly convex in $x$.

Let us also consider the active power loss function as a function of nodal power injections. By implicit function theorem using (1), state variables could be expressed through nodal injections:

$$F(x(S), S) \equiv 0,$$

and hence one obtains an implicit expression of the active power losses as a function of nodal injections:

$$L(S) \equiv L(x(S)).$$

However one can obtain explicit expressions for the derivatives of $L(S)$. By implicit function theorem the first order derivatives are given by the formula

$$\nabla_S L(S) = (J(x)^T)^{-1} \nabla_x L(x)$$

where $J(x)$ is the Jacobian of (1) with respect to $x$. The components of the gradient of $L(S)$ have a natural physical interpretation: the marginal increment $d$ of nodal load at bus $r$ will require the offsetting injection of $(1-\partial L(S)/\partial S_r)d$ at the swing bus to keep the system balanced. Here $-\partial L(S)/\partial S_r$ is the marginal loss coefficient (the sign depends on the convention what net injection is considered to be positive: here net generation is positive).

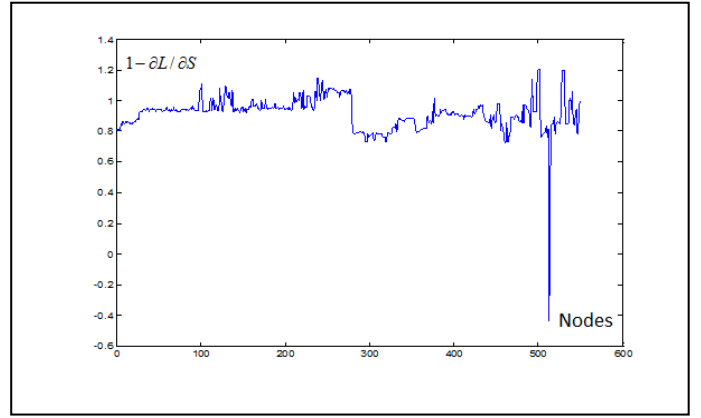

Fig. 1. The drawdown of the marginal loss coefficients indicate local excess of generation

The important property of $\nabla_S L(S)$ following from its very definition is that its components tend to infinity as $S$ approaches the solvability boundary and the Jacobian $J(x)$ tends to degeneracy. Use of this property of $\partial L(S)/\partial S_r$ for the analysis of power flow is also stressed in [7]. Moreover in many cases of practical interest where disturbances that bring the system towards the solvability boundary are local, the set of components of $\nabla_S L(S)$ that tend to infinity is also concentrated exactly at the area where the disturbances occur. Figure 1 illustrates this observation. The graph shows a fragment of a series of marginal loss coefficients for the 9000-bus model of the Russian grid. The extremely low levels of the coefficients (even below zero) in a number of closely located nodes indicate excess generation at this area. In this specific example the modelled line outage made the system unsolvable as the transmission capability was insufficient to provide output of a power plant.

The explicit expression for the hessian of $L(S)$ can also be obtained. Under natural assumptions $\nabla^2 L(S)$ is positive definite in some neighborhood of the origin.

## III. LOSS FUNCTION BASED REGULARIZATION

Consider the problem of minimization of the squared norm of the residual of the load flow equations regularized by the active power loss function:

$$\min_x 1/2 \| F(x,S) \|^2 + \alpha L(x) \quad (2)$$

where $\alpha > 0$ is a regularization parameter and $S$ is a given vector of nodal injections.

*Property 2.* If $x^*$ is a stationary point of (2) the residual of (1) at $x^*$ is proportional to the marginal loss coefficients:

$$F(x^*, S) = -\alpha (J(x^*)^T)^{-1} \nabla_x L(x^*) = -\nabla L(\overline{S}).$$

where $\overline{S}$ is the "corrected" nodal injections profile, such that

$$\bar{S} - S = F(x^*, S).$$

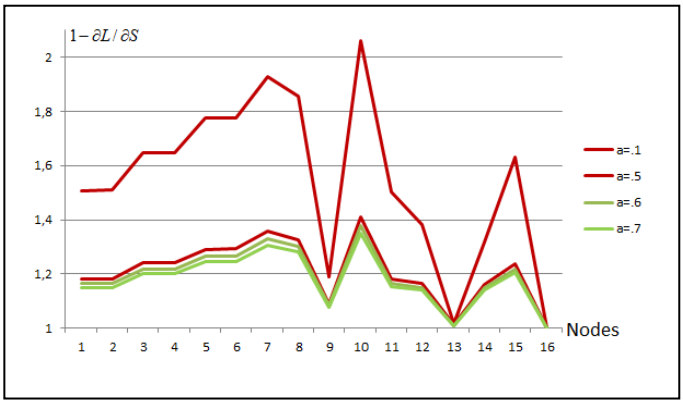

Fig. 2. Active power disturbance in 7 and 10 nodes: Marginal loss coefficients corresponding to active power injections

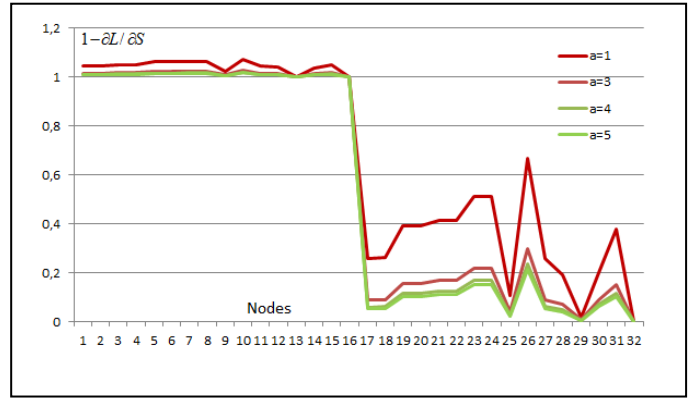

Fig. 3. Reactive power disturbance at node 10: Marginal loss coefficients corresponding to active and reactive power injections

As mentioned before the marginal loss coefficients while being relatively close to zero under normal conditions, tend to "mark" the problem areas by extremely high absolute values. Thus, Property 2 suggests that the sequence of iterates of a minimization method converges to a limit where the residual profile would indicate the problem location by the highest by absolute value nodal error.

However, the convergence critically depends upon convexity of the goal function in (2). It turns out that its convexity (at least locally at a stationary point) is closely related to the properties of the active losses as a function of nodal injections $L(S)$.

*Property 3.* If $E + \alpha \nabla^2 L(\bar{S})$ is positive definite (with $E$ being unity matrix) then a stationary point $x^*$ is a local minimum of (2).

As said previously $L(S)$ is convex at a neighborhood of the origin and hence if $\bar{S}$ happens to be at this neighborhood condition of Property 4 holds. To enforce convexity the regularization parameter $\alpha$ could be increased. The tradeoff, however, is the "quality" of identification of the failure location in the system. The following model examples illustrate the regularization.

## IV. EXAMPLE

We consider a model 16-bus grid with two types nodal load perturbations that make the load flow problem infeasible. The first type is excessive active power consumption at buses 7 and 10, and the second is excessive reactive power consumption at bus 10. The results of application of the loss function regularization for different parameter $\alpha$ are shown on Figure 1 and 2 for the corresponding types of perturbations. The graph on Figure 1 shows the profile of the marginal loss coefficients at the solution (which by Property 3 are proportional to the residual). It is clearly seen that maximum absolute values of the coefficients are attained at the buses where perturbations were made. However the profile tends to be smoother as the magnitude of the regularization parameter increases. All runs are made using the flat start as initial approximation. The red curves mark the cases where $\nabla^2 L(S)$ was indefinite at the solution. Despite the load correction that restores solvability of the load flow is obtained, the voltage magnitudes in the "red" cases are about 0.75 for some buses. On the contrary, the green curves correspond to the positive definite $\nabla^2 L(S)$ at the solution. The corresponding voltage levels are above 0.9.

Figure 2 shows that active power loss function regularization can also indicate locations with abnormal reactive power perturbations. One can see that the marginal loss coefficients w.r.t. reactive power injections take the highest values in the vicinity of the perturbed bus. As in Figure 1 green curves correspond to the positive definite $\nabla^2 L(S)$ at the solution and the red ones – to indefinite.

## V. CONCLUDING REMARKS

The experiments made so far with the loss function based regularization including the large scale grid model demonstrate potential efficiency of the approach in locating the source of the problem in the system. However there are still a number of research and practical questions. The best choice of the regularization parameter is one of them. The other is appropriate building of the regularization procedure into the market scheduling process in order to assure a "smooth" termination in case of infeasibility with the corresponding guidance for the operating personnel.

REFERENCES

[1] Davidson, M.R., Dogadushkina, Yu. V., Kreines, E. M., Novikova, N. M., Udaltsov, Yu. A., and Shiryaeva, L. V. "Mathematical Model of the Competitive Wholesale Power Market In Russia," J. of computer and Systems Sciences International; 43(3), pp. 72-83, 2004.

[2] Report on UPS of Russia operation 2014. System Operator of the United Power System, 2015. http://www.so-ups.ru/fileadmin/files/company/reports/disclosure/2015/ups_rep2014.pdf [in Russian]

[3] I.A. Hiskens and R.J. Davy. "Exploring the power flow solution space boundary," IEEE Trans. Power Systems, 16:389–95, 2001.

[4] D. Bienstock, "Progress on solving power flow problems," OPTIMA 93, pp. 1-7, 2013.

[5] A. Castillo and R.P. O'Neill, "Survey of approaches to the ACOPF," Optimal Power Flow Paper 4, 2013.

[6] T. J. Overbye. "Computation of a practical method to restore a power flow solvability," IEEE Trans. Power Systems, 10:280–286, 1995.

[7] N. Chemborisova. "Generalized performance indices in estimating power flow feasibility in deficit energy systems," Electricity, 5, pp. 2-10, 2004 [in Russian]